\documentclass[12pt]{article}
\usepackage{amssymb,amsbsy,amsmath,amsfonts,amssymb}
\usepackage{latexsym,euscript,exscale}

\usepackage{times}

\newcommand{\pf}{

\smallskip

\noindent {\it Proof : }}

\newcommand{\ext}{\preccurlyeq}

\newcommand {\kapa}{{\cal K}}

\newcommand {\N}{\mathbb N}

\newcommand {\R}{\mathbb R}

\newcommand{\pff}{$\hfill \square$
\smallskip}

\newcommand{\norm}[1]{\ensuremath{\left\|#1\right\|}}

\newtheorem{prop}{Proposition}[section]
\newtheorem{lemm}[prop]{Lemma}
\newtheorem{theo}[prop]{Theorem}

\newtheorem{defi}[prop]{Definition}
\newtheorem{conj}[prop]{Conjecture}

\title{On the number of non permutatively equivalent basic 
sequences in a Banach
space}

\author{Valentin Ferenczi}

\date{}

\begin{document}

\maketitle

\begin{abstract} 
 Let $X$ be a Banach space with a Schauder 
basis $(e_n)_{n \in \N}$. The relation $E_0$ is Borel
reducible to permutative equivalence between normalized 
block-sequences of $(e_n)_{n \in \N}$ or $X$ is
$c_0$ or $\ell_p$ saturated for some $1 \leq p<+\infty$. If
$(e_n)_{n \in \N}$ is shrinking unconditional then either it is equivalent to the
canonical basis of
$c_0$ or
$\ell_p, 1<p<+\infty$, or the relation
$E_0$ is Borel reducible to permutative
  equivalence between sequences of normalized disjoint blocks of $X$ 
or of $X^*$. If $(e_n)_{n \in \N}$ is unconditional, then either $X$ is isomorphic 
to $\ell_2$, or $X$
contains $2^{\omega}$ subspaces or $2^{\omega}$ quotients which are 
spanned by pairwise non
permutatively equivalent normalized unconditional bases.\footnote{MSC-class numbers: 46B03;
03E15}
\end{abstract}

\

\section{Introduction}

In the 1990's, W.T. Gowers and R.
Komorowski -
N. Tomczak- Jaegermann solved the so-called
Homogeneous Banach Space Problem. A Banach space is said to be
homogeneous if it is
isomorphic to its infinite dimensional closed subspaces; it is a 
consequence of two
theorems proved by these authors that a homogeneous Banach space
must  be isomorphic
to
$\ell_2$ \cite{G,KT}.

\

It is then natural to ask how many non isomorphic subspaces must 
contain a given Banach
space which is not isomorphic to $\ell_2$. This question was first asked the
author by G. Godefroy, and not much was known until 
recently about it in the literature, even concerning the classical spaces 
$c_0$ and $\ell_p$.

The correct
setting for this question is the classification of analytic 
equivalence relations on
Polish spaces by Borel reducibility. This area of research originated 
from the works of H.
Friedman and L. Stanley \cite{FS} and independently from the works of 
L. A. Harrington, A.
S. Kechris and A. Louveau \cite{HKL}. It may be thought of as an 
extension of the notion
of cardinality in terms of complexity, when one counts equivalence classes.

  If
$R$ (resp.
$S$) is an
equivalence relation on a Polish space $E$ (resp. $F$), then it is said
that $(E,R)$ is Borel
reducible to $(F,S)$ if there exists a Borel map $f:E \rightarrow F$ such that
$\forall x,y \in E, xRy \Leftrightarrow f(x)Sf(y)$.
An important equivalence relation is the relation $E_0$: it is
defined on  $2^{\omega}$ by
$$\alpha E_0 \beta \Leftrightarrow \exists m \in \N \forall n \geq m,
\alpha(n)=\beta(n).$$

The relation $E_0$ is a Borel equivalence relation with $2^{\omega}$ 
classes and which,
furthermore, admits
no Borel classification by real numbers, that is, there is no 
Borel map $f$
from
$2^\omega$ into $\R$ (equivalently, into a Polish space), such 
that $\alpha E_0 \beta
\Leftrightarrow f(\alpha)=f(\beta)$; such a relation is said to be 
{\em non-smooth}.
In fact $E_0$ is the $\leq_B$ minimum  non-smooth Borel equivalence 
relation \cite{HKL}.

  There is a natural way to equip the set of subspaces
of a Banach space
$X$ with a Borel structure (see i.e. \cite{Ke}), and the relation of isomorphism 
is analytic in this
setting \cite{B}. The relation
$E_0$ then appears as a natural threshold for results about isomorphism
between separable Banach spaces.
A Banach space $X$ was defined in \cite{FR2} to be {\em ergodic} if $E_0$ is Borel reducible
to isomorphism between subspaces of $X$; in particular, an ergodic Banach space
has continuum many non-isomorphic subspaces, and isomorphism between its
subspaces is non-smooth.

\

The question of the complexity of isomorphism between subspaces of a given Banach space $X$
is related to results and questions of Gowers about the structure of the relation of
embedding between subspaces of $X$
\cite{G}. In that article, Gowers proves the following structure theorem:

\begin{theo}(W.T. Gowers) Any Banach space contains a subspace $Y$ satisfying one of
the following properties, which are mutually exclusive and all possible:

\begin{itemize}

\item (a) $Y$ is hereditarily indecomposable (i.e. contains no direct sum of infinite
dimensional subspaces),
\item (b) $Y$ has an unconditional basis and no disjointly supported subspaces of $Y$ are
isomorphic,
\item (c) $Y$ has an unconditional basis and is strictly quasi-minimal (i.e. any two
subspaces of $Y$ have further isomorphic subspaces, but $Y$ contains no minimal subspace)
\item (d) $Y$ has an unconditional basis and is minimal (i.e. $Y$ embeds into any of its
subspaces).

\end{itemize}

\end{theo}

Note that these properties are preserved by passing to block-subspaces (in the associated
natural basis). Furthermore, knowing that a space belongs to one of the classes (a)-(d)
gives a lot of informations about operators and isomorphisms defined on it (see
\cite{G} about this).

\

C. Rosendal proved that any Banach space satisfying (a) is ergodic \cite{R}.
 The author and Rosendal noticed that a result of B. Bossard adapts easily to obtain that a
space satisfying (b) is ergodic \cite{FR2}. Finally by \cite{F}, using a result of \cite{R},
a space with (c) must be ergodic as well. 

It is furthermore known that a non-ergodic space $Y$ satisfying (d)  must be
isomorphic to its hyperplanes and to its square \cite{R}, must be reflexive, by \cite{FG}
and the classical theorem of James, and that it must
contain a block-subspace $Y_0$ such that
$Y_0
\simeq Y_0
\oplus Z$ for any block-subspace
$Z$ of
$Y$ \cite{FR2}.

Note that the class (d) contains the classical spaces $c_0$ and $\ell_p, 1 \leq p
<+\infty$, and Schlumprecht's space $S$ \cite{S}. Concerning those spaces, it is known that
$c_0$ and
$\ell_p$,
$1 \leq p <2$ are ergodic \cite{FG}. For $2<p<+\infty$, it is only known that there exist
$\omega_1$ non-isomorphic subspaces of $\ell_p$ (see \cite{LT}, Th. 2.d.9). The case of $S$
is also unsolved.

These
results suggest the following conjecture:

\begin{conj}\label{l2} Every separable Banach space is either
isomorphic to $\ell_2$ or ergodic.
\end{conj}

Now the spaces $c_0$
or $\ell_p$, $p \neq 2$ are also very homogeneous in some sense, since they are
isomorphic to any of their block-subspaces (with respect to their canonical basis).

 It also turns out that all the mentioned results about ergodic Banach
spaces (except of course \cite{FG}), as well as Gowers' theorem, can be proved
using block-subspaces of a given basis. So it is  natural to study the
homogeneity question restricted to block-subspaces of a Banach  space $X$ with a Schauder
basis. Block-subspaces can be thought of as "regular"  subspaces in this
context, for example, they will have a canonical unconditional 
basis, whenever the basis of
$X$ is unconditional.

In fact, classical results show that we can get a lot of information 
about the properties of
a space with a basis from the properties of its block-subspaces. For example, recall
that two basic sequences $(x_n)$ and $(y_n)$ are said to be equivalent
if the linear map $T$ defined on the closed linear span of $(x_n)$ by
$Tx_n=y_n, \forall n \in \N$ is an isomorphism onto the closed linear span
of $(y_n)$.   The canonical bases of $c_0$ and $\ell_p$
are characterized, up to equivalence of basis, by the property of being equivalent to all
their normalized block-bases (this  is Zippin's theorem,
\cite{LT} Th. 2.a.9).

\

If the basis is unconditional, it will also be natural to consider sequences of blocks (i.e. finitely supported
vectors) whose supports are disjoint, but not necessarily successive (equivalently,
block-sequences
 of  permutations of the basis). 
 This 
distinction is relevant
as some classical results require considering such basic sequences 
instead of 
block-sequences: for example,
\cite{LT} Theorem 2.10, according to which $c_0$ and $\ell_p$ are 
characterized by unconditionality and the property
that every subspace with a basis of disjointly supported blocks is 
complemented.

\

  We also note that
the Theorem of
Komorowski and Tomczak-Jaeger\-mann \cite{KT} is totally irrelevant in this context: it
shows the existence of an  "exotic" subspace of a
Banach space $X$ spanned by an unconditional basis, which has a unconditional
finite-dimensional decomposition but 
which fails to have an
unconditional basis, so it will give no information whatsoever on 
block-sequences or disjointly supported blocks of $X$.

\

The natural conjecture concerning the spaces $c_0$ and $\ell_p$ is as follows:

\begin{conj} if $X$ is a Banach space with an (unconditional)
  basis, is it true that either
   $X$ is isomorphic to its
block-subspaces or  $E_0$ is Borel reducible to isomorphism between the
block-subspaces of $X$? Is it true that if $X$ is isomorphic to its 
block-subspaces
then $X$ is isomorphic to $c_0$ or $\ell_p$? Are these assertions true 
when one replaces
block-subspaces by subspaces supported by disjointly supported blocks?
\end{conj}

Note
that by an easy result of \cite{FR1} using the theorem of Zippin, our
conjecture is solved if one replaces isomorphism by equivalence: if
$X$ is a Banach space with a normalized basis $(e_n)_{n \in \N}$, then either
$(e_n)_{n \in \N}$ is equivalent to the canonical basis of $c_0$ or
$\ell_p, 1 \leq p
<+\infty$, or $E_0$ is Borel reducible to equivalence between normalized
block-sequences of $X$.

  \

Some remarks and partial answers to these conjectures may be found in \cite{F2}.
As solving these questions seems to be out of reach for the moment, in this paper we
shall concentrate our efforts on the corresponding conjectures obtained by
replacing isomorphism by permutative equivalence. As it turns out, we
shall get results which are very close to  positive answers in that
case. Two basic sequences
$(x_n)_{n \in \N}$ and
$(y_n)_{n \in \N}$ are said to be  permutatively equivalent if
there is a permutation $\sigma$ on $\N$ such that
$(x_n)_{n \in \N}$ is equivalent to $(y_{\sigma(n)})_{n \in \N}$, in which case
we write $(x_n) \sim^{perm} (y_n)$.
Permutative equivalence
between Schauder bases is implied by equivalence and implies isomorphism of
the closed linear spans.

It is common to study permutative equivalence between normalized unconditional
basic sequences, since then any permutation of the basis
 is
again  a basic sequence. 
However some of our results will concern the general case
of permutative equivalence
between normalized basic sequences which are not necessarily unconditional.

\

We list several reasons for which studying permutative equivalence 
is relevant. First,
some classical results which are false or unknown for isomorphism can be proved for
permutative equivalence.
 The Theorem of Zippin admits a generalization to
permutative equivalence, due to Bourgain, Cazazza, Lindenstrauss
and Tzafriri:
if an unconditional basis is permutatively equivalent to all
its normalized block-sequences,
then it must be equivalent to the canonical basis of $c_0$ or $\ell_p$ 
\cite{BCLT}. Also, a Schroeder-Bernstein result is valid for
 permutative equivalence: if $(x_n)_{n \in \N}$ and $(y_n)_{n \in \N}$
 are unconditional
 basic
sequences such that each one is permutatively equivalent to a subsequence of
the other, then $(x_n)_{n \in \N}$ and $(y_n)_{n \in \N}$ are permutatively
equivalent (apparently first proved by Mityagin, \cite{M}, and
\cite{W,Wo}). Note that this is false without the unconditionality
asssumption, by the example of Gowers and Maurey of a space isomorphic to its
subspaces of codimension $2$, by a double shift of its natural basis, but not 
isomorphic to its hyperplanes \cite{GM2}. The Schroeder-Bernstein Problem for
Banach spaces, which asks whether two Banach spaces which are isomorphic to complemented
subspaces of each other must be isomorphic, is unsolved in the case of 
them having an unconditional basis, and solved by the negative in the general
case, by Gowers \cite{G2} and the examples of \cite{GM2}.

On the other hand,
permutative equivalence is already a complex relation. As isomorphism,
it is analytic non Borel, as we shall prove in Proposition 
\ref{bossard},
while equivalence of basic sequences is only $K_{\sigma}$ \cite{R1}. 
In fact, as far as we know, permutative equivalence between
basic sequences could well be as complex as isomorphism between Banach spaces
with a Schauder basis, or between separable Banach spaces in general.

 Also, some results of uniqueness 
of unconditional bases (see \cite{BCLT,K}, and \cite{CK1,CK2})
make it possible, in some special cases, to deduce permutative equivalence of 
basic sequences from
isomorphism of the Banach spaces they span. For example, the results of
\cite{FG} about the complexity of isomorphism, are essentially 
results about the complexity
of permutative equivalence: indeed, their constructions always 
realize a reducing of some
 equivalence relations to isomorphism between some subspaces equipped with 
canonical unconditional bases,
which are isomorphic exactly when these canonical bases are 
permutatively equivalent (\cite{FG} Theorem 2.6, Theorem 3.3).
The same holds in \cite{Rthese}, where it is used that subsequences of the basis of
Tsirelson's space are (permutatively) equivalent if and only if they span isomorphic
subspaces.

\

In this article, we investigate the complexity of permutative
equivalence between normalized basic sequences of a given Banach space; in particular, if a 
Schauder basis is not
equivalent to
$c_0$ or $\ell_p$, we ask how many non permutatively equivalent 
normalized block-sequences (resp. sequences
of disjointly supported blocks) it must contain.

\begin{conj}\label{conjectureperm} Let $X$ be a Banach space with a (resp. unconditional)
basis  which is not equivalent
to the canonical basis of $c_0$ or $\ell_p, 1 \leq p<+\infty$. Then
$E_0$ is Borel reducible to permutative equivalence between normalized
block-sequences (resp. sequences of disjointly supported blocks) of $X$.
\end{conj}

\

In Section 1, we extend the results of \cite{B} to prove that the relation of
permutative equivalence is non Borel, and the results of
\cite{FG} to show that it reduces the relation $E_{K_{\sigma}}$, and thus is not reducible
to the Borel action of a Polish group on a Polish space (Proposition \ref{bossard}).

 In Section 2, we prove several lemmas to obtain a result which is very
close to a
 positive answer to Conjecture \ref{conjectureperm}. If $X$ is a Banach space
 with a Schauder basis such that $E_0$ is not Borel reducible to
permutative
 equivalence between normalized block-sequences of $X$, then there exists
$p \in [1,+\infty]$ such that $X$ is $\ell_p$-saturated (or $c_0$-saturated if $p=+\infty$),
Theorem \ref{resultatblock}. If the basis is unconditional, then in fact
any normalized block-sequence of
$X$ has a
 subsequence which is equivalent to the canonical basis of $\ell_p$ (or $c_0$ if
$p=+\infty$), Theorem \ref{resultat}. If the basis is unconditional and
$E_0$ not Borel reducible to
permutative equivalence between normalized sequences of disjointly supported
 blocks,
then we also have that $p$ is unique such that $l_p$ is disjointly finitely representable on
$X$, and that $X$ satisfies an upper $p$ estimate (Theorem
 \ref{resultat}).

Our main tools for this result are a technical lemma (Lemma
 \ref{lppluslqtechnical}); a result of Rosendal about reducings of $E_0$ to
 equivalence relations between subsequences of a given basis (\cite{R},
 Proposition 21), which uses a result of Bourgain, Casazza, Lindenstrauss,
 Tzafriri \cite{BCLT};  Krivine's Theorem \cite{L} about finite block
 representability of the spaces $\ell_p$, and a result of stabilization of Lipschitz
 functions, by Odell, Rosenthal and Schlumprecht \cite{OR}.

In Section 3, we deduce that if  
 $X$ is a Banach space with a shrinking normalized unconditional  
 basis
$(e_n)$, then either $(e_n)$ is equivalent to the canonical basis of $c_0$ or 
some $\ell_p, 1<p<+\infty$,
or
$E_0$ is Borel reducible to permutative
equivalence between normalized disjointly supported sequences of blocks on $X$,
or on $X^*$ (Theorem \ref{dsb}).
It follows that if $X$ is a Banach space with an  unconditional basis, then
  either $X$  is  isomorphic
to $\ell_2$, or
$X$ contains
$2^{\omega}$ subspaces or $2^{\omega}$ quotients spanned by 
 unconditional bases which
are mutually non permutatively equivalent (Theorem \ref{continuum}).

\subsection{Notation}

Let us fix or recall some notation. For the reader interested in more details,
we refer to \cite{LT}.

A sequence $(e_n)_{n \in \N}$ with closed linear span $X$ is said to be ba\-sic
(or a Schau\-der basis of $X$)
if for any $x \in X$, there exists a unique scalar sequence $(\lambda_n)_{n
  \in \N}$ such that $x=\sum_{n \in \N}\lambda_n e_n$. 
This is equivalent to saying that there exists $C \geq 1$ such that for any
 $x=\sum_{n \in \N}\lambda_n e_n$, any integer $m$, 
$\norm{\sum_{n \leq m}\lambda_n e_n} \leq C \norm{x}$.
An interval of integers $E$ is the intersection of an  interval of $\R$ with
 $\N$; it will also denote the canonical projection on the span of $(e_n)_{n
   \in E}$, called interval projection. A Schauder basis is said to be
 bimonotone if every non-zero interval projection on its span is of norm $1$.
 A Banach space with a Schauder basis may always be renormed with an
equivalent norm so that the basis is bimonotone in the new norm.

Let $X$ be a Banach space with a Schauder basis $(e_n)_{n \in \N}$.
We shall  use some standard notation about blocks on $(e_n)_{n \in \N}$, i.e. finitely
supported non-zero vectors, for
example, we shall write $x<y$ and say that $x$ and $y$ are successive
when $\max(supp(x)) <
\min(supp(y))$.

The set of normalized block-sequences, i.e. infinite sequences of
 successive normalized blocks, in $X$ is denoted
$bb(X)$.
The set of normalized sequences of disjointly supported blocks in $X$ 
is denoted $dsb(X)$.
Both are seen here as metric spaces as subspaces of $X^\omega$ with the 
product of the norm
topology, and this turns them into Polish spaces.

If $(x_n)_{n \in I}$ is a finite or infinite sequence in $X$
then $[x_n]_{n \in I}$
will stand for its closed linear span.
  We recall that two basic sequences $(x_n)_{n \in \N}$ and
$(y_n)_{n \in \N}$ are said to be
{\em equivalent} if the map $T:[x_n]_{n \in \N} \rightarrow [y_n]_{n
\in \N}$ defined by
$T(x_n)=y_n$ for all $n \in \N$ is an isomorphism, in which case we write
$(x_n) \sim (y_n)$; if
$\norm{T}\norm{T^{-1}} \leq C$, then
  they are {\em $C$-equivalent}, and we write $(x_n) \sim^C (y_n)$.
  A basic sequence is
said to be {\em ($C$-)subsymmetric} if it is ($C$-)equivalent to all its
subsequences. Note that a subsymmetric sequence need not be unconditional. A Banach
space with
 a subsymmetric Schauder basis may always be
renormed to become $1$-subsymmetric. Two basic sequences
$(x_n)_{n \in \N}$ and
$(y_n)_{n \in \N}$ are said to be {\em permutatively equivalent} if
there is a permutation $\sigma$ on $\N$ such that
$(x_n)_{n \in \N}$ is equivalent to $(y_{\sigma(n)})_{n \in \N}$, in which case
we write $(x_n) \sim^{perm} (y_n)$.

 Let
$c_{00}$ denote the set of eventually null scalar sequences.
If $(x_n)_{n \in I}$ and $(y_n)_{n \in I}$ are finite or infinite 
 basic sequences,
we shall say that {\em $(y_n)$ $C$-dominates $(x_n)$}, and write
$(x_n) \leq^C (y_n)$, to mean that for all $(\lambda_i)_{i \in I}$ in
$c_{00}$, $\norm{\sum_{i \in I}\lambda_i x_i} \leq C
 \norm{\sum_{i \in I}\lambda_i y_i}$.

A basic sequence $(u_i)_{i \in \N}$ is said to be $C$-unconditional if for any sequence
of signs $(\epsilon_i)_{i \in \N} \in \{-1,1\}^{\omega}$, any sequence 
$(\lambda_i)_{i \in \N} \in
c_{00}$, we have
$\norm{\sum_{i \in \N} \epsilon_i \lambda_i u_i} \leq C \norm{\sum_{i \in \N}
\lambda_i u_i}$. In particular, any canonical projection on the 
closed linear span of some
subsequence of a $1$-unconditional basis is of norm $1$.  We may always assume by 
renorming that an
unconditional basis is $1$-unconditional.  If in addition the basis is subsymmetric, we may  ensure
that  it is also $1$-subsymmetric
in the new norm.

\

\subsection{General results about permutative equivalence}

In this introductive section, we recall the setting defined by B. Bossard for 
studying the complexity of equivalence relations between basic sequences, and
 notice that his results about 
isomorphism easily extend to
permutative equivalence \cite{B}.

Let $u$ be the normalized universal basic sequence of Pe\l czy\'nski \cite{P} and $U$ be
its closed linear span. The sequence $u$ is defined by  the
following property: any normalized basic sequence in Banach space is equivalent to a
subsequence $u'$ of $u$ such that the canonical projection from $U$ onto the
span of $u'$ is bounded.

Bossard defined a natural coding of basic sequences by considering the subsequences of $u$ (identified with infinite subsets of
$\N$).   Thus a property of basic sequences is Borel (resp. analytic,...) if the set of
subsequences of $\N$, canonically identified with subsequences of $u$,
 with this property is a Borel (resp. analytic...) subset of $[\omega]^{\omega}$ (the set of
increasing sequences of integers).

 The sequence $u$ also has an unconditional version $v=(v_n)_{n \in \N}$,  i.e. $v$ is
a normalized
 unconditional basic sequence and any normalized
 unconditional basic sequence in a Banach space is equivalent to a subsequence
 of $v$. We may represent $v$ as a subsequence of $u$.

\

The relation $E_{K_{\sigma}}$ is defined as the maximum $K_{\sigma}$ relation
on a Polish space
for the order $\leq_B$ of Borel reducibility \cite{R1}.  For details in the Banach space
context we refer to
\cite{FG}; let us just note here that $E_{K_{\sigma}}$ can not
(and thus neither can a  relation which
Borel reduces it)  be reduced to the Borel action of a Polish
group on a Polish space.

\begin{prop}\label{bossard} The relation of permutative equivalence 
between normalized basic
sequences is analytic non Borel and it Borel reduces $E_{K_{\sigma}}$.
In particular it cannot be Borel reducible to a relation associated to
the Borel action of a Polish group on a Polish space. 
\end{prop}

\pf By \cite{FG}, the relation $E_{K_{\sigma}}$ is Borel reducible
to isomorphism between Banach spaces. In the list of equivalence of
\cite{FG} Theorem 2.6, we may obviously add the condition: "is permutatively
equivalent to", since equivalence of bases implies permutative equivalence
which in turn implies isomorphism of the closed linear spans. This implies that
$E_{K_{\sigma}}$ is Borel reducible to permutative equivalence.
Note that the reduction of $E_{K_{\sigma}}$ is obtained using
unconditional 
sequences in
$\ell_p, 1 \leq p<2$ (resp. $c_0$), and so $E_0$ is Borel reducible to permutative
equivalence between some unconditional sequences
 in $\ell_p, 1 \leq p<2$ (resp. $c_0$),
 and in particular $\ell_p, 1 \leq p<2$ (resp. $c_0$) contains $2^{\omega}$
non permutatively equivalent unconditional basic sequences. 
This fact will be used at the end of this article.

It is immediate that permutative equivalence is analytic (this was already
observed in \cite{FR1}).
To prove that it is not Borel, we now define an unconditional version of a family of
basic sequences indexed by the set $\cal T$ of trees on $\omega$, which was considered
in \cite{B}. We also refer to \cite{B} for more details about the proof or the notation, in
particular concerning trees.

Let $T=\omega^{<\omega}$ denote the set of finite sequences of integers. Let $c_{00}(T)$ be
the space of finitely supported functions from
$T$ to
$\R$ and let $\phi_s:T \rightarrow \{0,1\}$ be the characteristic function of $\{s\}$ for
every $s \in T$. An admissible choice of intervals is a finite set $\{I_j, 0 \leq j \leq
k\}$ of intervals of $T$ such that every branch of $T$ meets at most one of these
intervals. We consider the $\ell_2$-James tree space $\tilde{v}(T)$ on $v$, i.e. the
completion of $c_{00}(T)$ under the norm defined by
$$\norm{y}=\sup((\sum_{j=0}^k\norm{\sum_{s \in I_j} y(s)v_{|s|}}^2)^{1/2}),$$
where $|s|$ is the length of $s \in T$ and the sup is taken over $k \in \N$ and all
admissible choices of intervals $\{I_j, 0 \leq j \leq k\}$.

If $A \subset T$, we let $\tilde{v}(A)$ be the subspace of $\tilde{v}(T)$ generated by
$\{\phi_s, s \in A\}$. We thus have defined a map $\tilde{v}$ on $\cal T$ to subsequences of
$v$ and thus of $u$. We claim that $\tilde{v}$ satisfies the following properties:

\begin{itemize}

\item a) $\tilde{v}$ is Borel, 

\item b) for all $\theta$, $\tilde{v}(\theta)$ is unconditional,

\item c) if $\theta$ is well-founded 
then $\tilde{v}(\theta)$ spans
a reflexive space,

\item d) if $\theta$ is ill-founded then some subsequence 
of $\tilde{v}(\theta)$
(corresponding to a branch of $\theta$) is equivalent to $v$.
\end{itemize}

The facts a), c) and d)
are valid for an $\ell_2$-James space on any Schauder basis instead of $(v_n)$. The proof
of a) is essentially the same as \cite{B} Lemma 2.4. Reproduce \cite{B} Lemma 1.5 and the
Fact in the proof of \cite{B} Theorem 1.2 for c), and \cite{B} Lemma 1.4 for d).

To prove b), we write an unconditional version of \cite{B} Lemma 1.3. Consider a real
sequence
$(\lambda_i)_{i
\in
\N}$,
$I$ an interval of
$T$, an integer $n \in \N$ and a subset $J$ of $[0,n]$. We denote by $c$ an upper bound for
the norms of canonical projections on subsequences of $v$. As in \cite{B}, $\kapa: \omega
\rightarrow \omega^{<\omega}$ is a fixed enumeration of $\omega^{<\omega}$ such that if $s
\ext s'$ then $\overline{s} \leq \overline{s'}$, if $\overline{s}$ denotes
$\kapa^{-1}(s)$. Write
$s_n=\kapa(n)$.

For $s \in T$, $(\sum_{i \in J} \lambda_i \phi_{s_i})(s)$ 
is equal to $\lambda_{\overline{s}}$ if $\overline{s} \in J$ and to $0$ otherwise.
Therefore,
$$\norm{\sum_{s \in I}(\sum_{i \in J}\lambda_i \phi_{s_i})(s)u_{|s|}}=
\norm{\sum_{s \in I, \overline{s} \in J}\lambda_{\overline{s}}u_{|s|}}$$
$$\leq c\norm{\sum_{s \in I, \overline{s} \leq n}\lambda_{\overline{s}}u_{|s|}}
=\norm{\sum_{s \in I}(\sum_{i \leq n}\lambda_i \phi_{s_i})(s)u_{|s|}},$$
since if $s,s' \in I$ then $s \ext s'$ iff $\overline{s} \leq \overline{s'}$.
Let  $\{I_j, 0 \leq j \leq k\}$ be an admissible choice of intervals. We have
$$\sum_{j=0}^k \norm{\sum_{s \in I_j}(\sum_{i \in J}\lambda_i \phi_{s_i})(s)u_{|s|}}^2 \leq
c^2 
\sum_{j=0}^k \norm{\sum_{s \in I_j}(\sum_{i \leq n}\lambda_i \phi_{s_i})(s)u_{|s|}}^2.$$
Thus $$\norm{\sum_{i \in J}\lambda_i \phi_{s_i}} \leq c\norm{\sum_{i \leq n}\lambda_i
\phi_{s_i}},$$
and $(\phi_{s_i})_{i \in \omega}$ is an unconditional basic sequence. The fact b) follows.

We note the following fact about $v$. If $v$ is  equivalent to the
subsequence of some normalized unconditional
basic sequence $w$, then $v$ is permutatively equivalent
to
$w$; indeed $w$ is equivalent to a subsequence of  $v$
 by definition of $v$ and the
result follows by the 
Schroeder-Bernstein's principle for permutative equivalence mentioned in the
 introduction \cite{M,W,Wo}. So it follows from b) and d):

\begin{itemize}
\item d') if $\theta$ is ill-founded then $\tilde{v}(\theta)$ is permutatively equivalent to
$v$.
\end{itemize}

 By c), $v(\theta)$ and $v$ are never permutatively equivalent when 
$\theta$ is well-founded. If $A$ is the $\sim^{perm}$-class of 
 $v$, it follows from this and from d') that
${\cal T}\setminus WF=v^{-1}(A)$, where $WF$ denotes the set of ill-founded trees on
$\omega$. So by a) and the well-known fact that $WF$ is non Borel, $A$ is non Borel, and it
follows that $\sim^{perm}$ is non Borel. \pff

We note here that the relations $=^+$, and the product $E_{K_{\sigma}} \otimes =^+$, defined
as in
\cite{FG}, may, by similar observations as in the $E_{K_{\sigma}}$ case, be reduced to
permutative equivalence between basic sequences.

\section{Reducing $E_0$ to permutative equivalence.}

\subsection{Reducing $E_0$ to permutative equivalence between block-sequences.}

Our initial and important technical
result bares similarity with \cite{LT} Lemma 2.a.11: from 
an hypothesis on
block-sequences of a Banach space, we already get a lot of 
information by looking at those
block-sequences of the form $((1-\lambda_n)x_n+\lambda_n y_n)_{n \in 
\N}$, for some fixed
sequences
$(x_n)$ and $(y_n)$ and choices of  sequences $(\lambda_n)_{n \in \N} 
\in [0,1]^{\N}$.

  Let
$(x_n)_{n \in \N}$ and
$(y_n)_{n \in \N}$ be normalized basic sequences generating spaces $X$ and $Y$. We equip
$X
\oplus Y$ with its canonical normalized basis
$(e_n)_{n \in \N}$,
that is, for any $(\mu_n)_{n \in \N} \in c_{00}$,
$$\norm{\sum_{n \in \N} \mu_n e_n}=\norm{\sum_{n \in \N} 
\mu_{2n-1}x_n}+\norm{\sum_{n \in
\N}
\mu_{2n}y_n}.$$

We shall identify vectors in $X$ (resp. $Y$) with their image in $X \oplus Y$.
Given a sequence $(a_n)_{n \in \N}
\in [0,1]^{\N}$, the sequence
  $(a_i x_{i}+(1-a_i)
y_i)_{i \in \N}$ is a normalized block-sequence of $X \oplus Y$.
We denote by $bb_{2}(X \oplus Y)$ the set of such infinite block-sequences.

Let $(I_k)_{k \in \N}$ be  a
sequence of successive intervals of integers forming a partition of $\N$, i.e.
$\forall k \in \N, \min I_{k+1}=\max I_k +1$, and let
$(\delta_k)_{k \in \N}$ be a positive decreasing sequence
  converging to $0$. We shall say that $(I_k),(\delta_k)$ is a 
{\em rapidly converging
system} if
$\delta_1 \leq 1/2$ and for all $k \geq 1$:

\begin{itemize}
\item (1) $|I_{k}|\delta_{k+1} \leq 1/4,$

\item (2) $|I_k|/2 >  \sum_{j=1}^{k-1}|I_j|.$

\end{itemize}

For any $\alpha \in 2^{\omega}$, we define a sequence of positive numbers
$(a_n(\alpha))_{n \in \N}$ by
$$a_n(\alpha)=\delta_{k+\alpha(k)},\ \forall k \in \N, \forall n \in I_k.$$

Finally we define a map $f$ from $2^{\omega}$ into
$bb_2(X \oplus Y)$ by 
$$f(\alpha)=(a_i(\alpha)x_i+(1-a_i(\alpha))y_i)_{i \in \N}.$$
We shall say that $f$ is the {\em map associated to} $(I_k),(\delta_k)$.

\begin{lemm}\label{lppluslqtechnical} Assume $X$ (resp. $Y$) is a 
Banach space with a
normalized Schauder basis
$(x_n)$ (resp. $(y_n)$). Let $(I_k)$, $(\delta_k)$ form a rapidly 
converging system  and
$f:2^{\omega} \rightarrow bb_2(X \oplus Y)$ be the associated map. Then
$f$ Borel reduces the relation
$E_0$ to permutative equivalence on $bb_2(X \oplus Y)$ or there exist 
$C \geq 1$, an
infinite subset $K$ of $\N$, and for each $k \in K$, a subset $J_k$ 
of $I_k$ with $|I_k
\setminus J_k|
\leq
\sum_{j=0}^{k-1}|I_j|$, and distinct integers $(n_i)_{i
\in J_k}$  such that
$$(\delta_k x_i+ y_i)_{i \in J_k}
\sim^C (y_{n_i})_{i \in J_k}.$$
\end{lemm}

\pf Without loss of generality we assume that
$(x_n)$ and $(y_n)$ are bimonotone.

The map $f$ is obviously Borel (even continuous) and whenever $\alpha 
E_0 \beta$,
$f(\alpha)$ is equivalent, and thus permutatively equivalent to $f(\beta)$.

Assume $f$ does not Borel reduce $E_0$ to permutative equivalence on 
$bb_2(X \oplus Y)$.
We have $f(\alpha) \sim^{perm} f(\beta)$ for some $\alpha,\beta$ in 
$2^{\omega}$ which
are not $E_0$ related, and let $C$ be the associated constant of equivalence.
We may assume for arbitrarily large $k$ that $\alpha(k)=0$ while 
$\beta(k)=1$. Let $K$ be
the infinite set of such integers, and let $k \in K$.

By the permutative equivalence between $f(\alpha)$ and $f(\beta)$, the
sequence $(\delta_k x_i+(1-\delta_k) y_i)_{i \in I_k}$ satisfies
$$(\delta_k x_i+(1-\delta_k) y_i)_{i \in I_k} \sim^C (\delta_{k_i} 
x_{n_i}+(1-\delta_{k_i})
y_{n_i})_{i \in I_k},$$
where
$(n_i)_{i \in I_k}$ is a sequence of distinct integers, and $\forall 
i \in J_k$, $k_i$ is
equal to
$m+\beta(m)$ if $m$ is such that $n_i \in I_m$.

By the increasing condition (2), there exists a subset $J_k$ of 
$I_k$, of length at least
$|I_k|-\sum_{j=1}^{k-1}|I_j|>0$,
for which we have
$$(\delta_k x_i+(1-\delta_k) y_i)_{i \in J_k} \sim^C (\delta_{k_i} 
x_{n_i}+(1-\delta_{k_i})
y_{n_i})_{i \in J_k},$$

where for $i \in J_k$, $k_i$ is of the form $m+\beta(m)$ for some $m 
\geq k$. Since $\beta(k)=1$ it follows 
that for all $i \in
J_k$,
$k_i
\geq k+1$ and thus
$\delta_{k_i} \leq \delta_{k+1}$.

By 
the previous remark, for any $(\lambda_i)_{i \in J_k}$,

$$\norm{\sum_{J_k}\delta_{k_i}\lambda_i x_{n_i}} \leq \delta_{k+1}|J_k| \max_{i \in 
J_k}|\lambda_i|,$$

so  as
$\delta_{k+1}|J_k|
\leq 1/4$, and by bimonotonicity,

$$\norm{\sum_{J_k}\delta_{k_i}\lambda_i x_{n_i}} \leq 1/4
\norm{\sum_{J_k}\lambda_i
y_{n_i}}.$$

By the same type of estimate, we have that

$$\frac{3}{4}\norm{\sum_{J_k}\lambda_i
y_{n_i}} \leq \norm{\sum_{J_k}(1-\delta_{k_i})\lambda_i
y_{n_i}} \leq \frac{5}{4}\norm{\sum_{J_k}\lambda_i
y_{n_i}}.$$

Finally, $(\delta_{k_i} x_{n_i}+(1-\delta_{k_i})y_{n_i})_{i \in J_k} 
\sim^3 (y_{n_i})_{i
\in J_k}$. Also,
$$\frac{1}{2}\norm{\sum_{J_k}\lambda_i
(\delta_k x_i+y_i)} \leq \norm{\sum_{J_k}\lambda_i
(\delta_k x_i+(1-\delta_k)y_i)} \leq \frac{3}{2}\norm{\sum_{J_k}\lambda_i
(\delta_k x_i+y_i)},$$

since $\delta_k \leq 1/2$, so $(\delta_{k} x_{i}+(1-\delta_{k})y_{i})_{i \in J_k} \sim^3
(\delta_k  x_i+y_i)_{i \in J_k}$,
and it follows

$$(\delta_k x_i+y_i)_{i \in J_k} \sim^{9C} (y_{n_i})_{i \in J_k}. \ 
\hfill \square$$

\

Let $\ext$ be a linear order  on $\N$. When $I$ is a finite 
subset of $\N$, we
denote by
$(I)_i^{\ext}$ the
$i$-th element of
$I$ written in $\ext$-increasing order.

\begin{defi} Let $(y_n)$ be a $1$-subsymmetric $1$-unconditional 
 basic sequence.
Let
$\ext$ be a linear order  on $\N$.  We define a normed space with a $1$-unconditional basis
$(y_n)^{\ext}$ by letting, for all $k \in \N$, for all 
$(\lambda_i)_{i=1}^k \in \R^k$,
$$\norm{\sum_{i=1}^k \lambda_i y_i^{\ext}}_{\ext}=\norm{\sum_{i=1}^k
\lambda_i y_{\{1,\ldots,k\}_i^{\ext}}}.$$
\end{defi}

\

We note a few easy facts. If $\leq$ is the usual order
relation on
$\N$, then $(y_n)^{\ext}$ is obviously
$1$-equivalent to
$(y_n)$. When $(y_n)$ is $1$-symmetric (i.e. $1$-equivalent to $(y_{\sigma(n)})$
for any permutation $\sigma$ on $\N$), then the sequence 
space defined by
$\norm{.}_{\ext}$ is always $1$-equivalent to $(y_n)$.  We shall also
be interested in
$\norm{.}_{\geq}$, where $\geq$ is defined as usual on $\N$; note 
that this defines a
$1$-subsymmetric basic sequence, and that $(y_n)^{\geq\geq}$ is 
$1$-equivalent to $(y_n)$.
We also note that the operation
$\ext$ preserves domination.

If $(y_n)$ is a subsymmetric unconditional basis, then we define 
$(y_n)^{\ext}$ as
$(y'_n)^{\ext}$, if $(y'_n)$ is the canonical $1$-subsymmetric 
$1$-unconditional basis
equivalent to $(y_n)$. The previous observations are still valid up 
to some constant of
equivalence.

\begin{prop}\label{lppluslq} Let $X$ be a Banach space with a 
normalized unconditional basis
$(x_n)$ and $Y$ be a Banach space with a normalized subsymmetric 
unconditional basis
$(y_n)$. The relation $E_0$ is Borel reducible to permutative
equivalence on $bb_2(X \oplus Y)$ or there exists a linear order $\ext$
on $\N$ such that $(y_n) \leq (y_n)^{\ext}$ and $(x_n) \leq (y_n)^{\ext}$.
\end{prop}

\pf Without loss of generality we assume that
$(x_n)$ is
bimonotone  and that
$(y_n)$ is $1$-unconditional and
$1$-subsymmetric. We
consider the following:

{\em Fact}: there exists $C \geq 1$ such that for all $n \in \N$, there
exists a permutation $\sigma_n$ of $\{1,\ldots,n\}$ such that
$(x_i)_{i=1}^n
\leq^C (y_{\sigma_n(i)})_{i=1}^n$ and $(y_i)_{i=1}^n \leq^C 
(y_{\sigma_n(i)})_{i=1}^n$.

\

We first assume the fact holds. For any $k \leq n$ we may define a linear order
$\ext^k_n$ on $\{1,\ldots,k\}$ by $i \ext_n^k \Leftrightarrow \sigma_n(i) \leq \sigma_n(j)$.

By the pidgeonhole principle, we may find for each $k$ some infinite set $M_k$ such that
forall $n \in M_k$, $\ext_n^{k}=\ext^k$ for some fixed linear order $\ext^k$ on
$\{1,\ldots,k\}$, and we may take care that $M_k \subset 
M_{k-1}$ for all $k$.

Therefore, whenever $i,j \leq n$, $i \ext^n j$ if and only if for some (equivalently for
all) $m \in M_n$, $\sigma_m(i) \leq \sigma_m(j)$. It follows that whenever $i,j \leq k \leq
n$, $i \ext^n j$ iff $i \ext^k j$.
We may therefore define a linear order $\ext$ on $\N$ by $i \ext j$ if and only if
$i \ext^n j$ for some (equivalently for all) $n \geq \max(i,j)$.

Since for any $k \in \N$ and $n \in M_k$, $(x_i)_{i=1}^k \leq^C
(y_{\sigma_n(i)})_{i=1}^k$, we conclude that
$(x_i)_{i=1}^k \leq^C (y_i^{\ext})_{i=1}^k$.
It follows that $(x_n) \leq^C (y_n)^{\ext}$. Likewise we
  obtain $(y_n) \leq^C (y_n)^{\ext}$.

\

Assume now the Fact does not hold.  We may build by induction a rapidly converging system
$(\delta_k),(I_k)$, so $\delta_1 \leq 1/2$ and for all $k \geq 1$:

\begin{itemize}
\item (1) $|I_{k}|\delta_{k+1} \leq 1/4,$

\item (2) $|I_k|/2 >  \sum_{j=1}^{k-1}|I_j|,$

\end{itemize}

and an increasing sequence of integers $(K_k)$ so that for all $k \geq 1$,

\begin{itemize}

\item (3) $K_k>\sum_{j=1}^{k-1} |I_j|$ and $K_k \delta_k \geq k$,

\item (4) for any permutation $\sigma$ on $\{1,\ldots,\max(I_k)\}$, 
there exists a sequence
$(\mu_i)_{i \leq \max(I_k)}$ of non-negative numbers with
$\norm{\sum_{i \leq \max(I_k)} \mu_i y_{\sigma(i)}} \leq 1$ and $\norm{\sum_{i
\leq \max(I_k)}
\mu_i x_i}+\norm{\sum_{i
\leq \max(I_k)}
\mu_i y_i}
\geq 5K_k$.

\end{itemize}

We note that all $\mu_i$'s in (4) are smaller than $1$. Also, any permutation on $I_k$
may be extended to a permutation on $\{1,\ldots,\max(I_k)\}$.
Thus using (3) and the bimonotonicity of the basis, we deduce from (4):

\begin{itemize}

\item (5) for any permutation $\tau$ on $I_k$, there exists a sequence
$(\mu_i)_{i \in I_k}$ of non-negative numbers such that $\norm{\sum_{i
\in I_k}
\mu_i x_i}+\norm{\sum_{i
\in I_k}
\mu_i y_i} 
\geq 3K_k$ and such that
$\norm{\sum_{i \in I_k} \mu_i y_{\tau(i)}} \leq 1$.

\end{itemize}

Now we claim that
the map associated to the system $(\delta_k),(I_k)$ Borel reduces
$E_0$ to permutative equivalence on $bb_2(X \oplus Y)$. Otherwise, 
by Lemma \ref{lppluslqtechnical}, we find $C \geq 1$, an
infinite subset $K$ of $\N$, and for all $k \in K$, a subset $J_k$ of 
$I_k$ with $|I_k
\setminus J_k|
\leq
\sum_{j=0}^{k-1}|I_j|$, and distinct integers $(n_i)_{i \in
J_k}$  such that, for any
$(\lambda_i)_{i \in J_k}$,
$$\delta_k (\norm{\sum_{i \in J_k}\lambda_i x_i}+\norm{\sum_{i \in 
J_k}\lambda_i y_i}) \leq
\delta_k \norm{\sum_{i \in J_k}\lambda_i x_i}+\norm{\sum_{i \in 
J_k}\lambda_i y_i}
\leq C\norm{\sum_{i \in J_k}\lambda_i y_{n_i}}.$$
Now by $1$-subsymmetry of $(y_n)$, the sequence $(y_{n_i})_{i \in 
J_k}$ is $1$-equivalent
to some $(y_{\sigma(i)})_{i \in J_k}$ for some permutation $\sigma$ of
$J_k$. We may extend $\sigma$ to a permutation $\tilde{\sigma}$ of
$I_k$.

   Applying the previous inequality to the coefficients
$\mu_i$ given by (5) for $\tau=\tilde{\sigma}$, we obtain

$$\delta_k (3K_k-2|I_k \setminus J_k|) \leq C \norm{\sum_{i \in J_k} \mu_i 
y_{\sigma(i)}},$$

so, by choice of $J_k$ and by $1$-unconditionality,

$$\delta_k (3K_k-2\sum_{j=1}^{k-1}|I_j|) \leq C \norm{\sum_{i \in I_k} \mu_i
y_{\tilde{\sigma}(i)}},$$

so by (3),
$$k \leq K_k \delta_k \leq C,$$ for arbitrary large $k \in K$, a contradiction.
\pff

\

\begin{prop}\label{b2} Let $X$ (resp. $Y$) be a Banach space with a normalized
subsymmetric unconditional basis
$(x_n)_{n \in \N}$ (resp. $(y_n)_{n \in \N}$). Assume
$(x_n)$ and $(y_n)$ are not equivalent. Then
 $E_0$ is  
Borel reducible to
permutative equivalence on
$bb_2(X
\oplus Y)$. \end{prop}

\pf Assume $(x_n)$ and $(y_n)$ are $1$-subsymmetric. We
assume $E_0$ is not Borel reducible to permutative equivalence on
$bb_2(X \oplus Y)$ and apply
Proposition
\ref{lppluslq}: let
$\ext$ be a linear order on
$\N$ such that
$(x_n)
\leq (y_n)_{\ext}$ and $(y_n) \leq (y_n)_{\ext}$.
  By a standard application of Ramsey's Theorem
for sequences of length
$2$, we may find an infinite subset
$K$ of
$\N$ on which either $\ext$ coincides with $\leq$ or $\ext$ coincides 
with $\geq$.

In the first case, by passing to a subsequence with indices in $K$, 
and by subsymmetry of
$(x_n)$ and $(y_n)$, we obtain that $(x_n) \leq (y_n)$.

In the second case, we have
$(x_n) \leq (y_n)^{\geq}$ and $(y_n) \leq (y_n)^{\geq}$. But this means that
$(y_n)^{\geq} \leq (y_n)^{\geq\geq}$, and as $(y_n)^{\geq\geq}$ is 
equivalent to
$(y_n)$, that
$(y_n) \geq (y_n)^{\geq}$. We deduce in that case that $(x_n) \leq (y_n)$ as well.

  By symmetry we obtain that these
two sequences are equivalent.
\pff

\

An immediate consequence of this fact is that $E_0$ is Borel reducible to
permutative equivalence between normalized block-sequences of $\ell_p \oplus
\ell_q$, $1 \leq p<q<+\infty$, and of $c_0 \oplus \ell_p$, $1 \leq p<+\infty$.

\

We recall a conjecture by H. P. Rosenthal. A Schauder basis
  $(e_n)_{n \in \N}$ is said to be a {\em Rosenthal basis} if any 
normalized block-sequence
of
$(e_n)_{n \in \N}$ has a subsequence which is equivalent to $(e_n)_{n 
\in \N}$. A Banach
space has {\em Rosenthal property} if it admits a Rosenthal basis.

It is not difficult to see that a Rosenthal basis must be 
subsymmetric unconditional. Also,
all spreading models generated by block-sequences are equivalent in a 
Banach space with a
Rosenthal basis. Rosenthal conjectured that any Rosenthal basis must 
be equivalent to the
canonical basis of
$c_0$ or $\ell_p$, $1 \leq p <+\infty$. For more details about this property, see 
\cite{FPR}.

\begin{lemm} \label{rosenthal} Let $X$ be a Banach space with an unconditional
 basis
$(e_n)_{n
\in \N}$. Assume $E_0$ is not Borel reducible to permutative
equivalence on $bb(X)$. Then there is a subsequence
$(f_n)_{n \in \N}$ of $(e_n)_{n \in \N}$ such that
every normalized block-sequence in $X$ has a subsequence which is
equivalent to $(f_n)_{n \in \N}$. In particular $(f_n)_{n \in \N}$ is 
a Rosenthal basic
sequence.
\end{lemm}

\pf Assume $E_0$ is not Borel reducible to permutative
equivalence on $bb(X)$. Then $E_0$ is Borel reducible to permutative
equivalence on the set of subsequences of $(x_n)_{n \in \N}$ for no 
$(x_n)_{n \in \N}$ in
$bb(X)$. By \cite{R} Proposition 21, it follows that every normalized block-sequence of
$X$ has a subsymmetric subsequence. It remains to show that 
any two subsymmetric
 block-sequences
$(x_n)$ and $(y_n)$ in
$X$ are equivalent. We may assume, by
passing to subsequences, that $x_k < y_k < x_{k+1}$ for all $k \in 
\N$. We then apply
Proposition \ref{b2}, since $E_0$ cannot be reduced to $\sim^{perm}$ on
$bb_2([x_k]_{k
\in
\N}
\oplus [y_k]_{k \in \N})$.
\pff

\

Let $X$ have a Schauder basis $(e_n)_{n \in \N}$.
For $1 \leq p \leq +\infty$, we say that {\em $\ell_p$ is block-finitely 
representable in $X$}
if there exists $C \geq 1$ such that $\forall n \in \N$, some length 
$n$ block-sequence in
$X$ is $C$-equivalent to the canonical basis of $\ell_p^n$. Note that 
this differs slightly
from the usual definition where it is required that we may take 
$C=1+\epsilon$ for any
$\epsilon>0$. By Krivine's theorem \cite{L}, there always exists $p \in 
[1,+\infty]$ such that
$\ell_p$ is block-finitely representable in $X$ (with $C$ arbitrarily close to
$1$ if you wish).
We say  that {\em $\ell_p$ is disjointly finitely 
representable in $X$}
if there exists $C \geq 1$ such that $\forall n \in \N$, some length 
$n$ sequence of disjointly supported blocks in
$X$ is $C$-equivalent to the canonical basis of $\ell_p^n$.

Using the proof by Lemberg of Krivine's Theorem \cite{L}, Odell, Rosenthal and
Schlum\-precht
\cite{OR} proved that  if $X$ is a Banach
space with a Schauder basis, $\oplus_{n \in \N} F_n$ is a decomposition of $X$ in
successive finite-dimensional subspaces of increasing dimension (where each $F_n$ is
equipped with the canonical basis  which is a subsequence of the basis of $X$),
$(\epsilon_n)$ is a sequence of positive reals, and $f: X 
\rightarrow \R$ is
a Lipschitz function on $X$, then there exists a subsequence $F_{k_n}$ of
$F_n$, finite block-subspaces $G_n$ of $F_{k_n}$ of increasing dimension,  and a map
$\tilde{f}$ on $\R^{<\omega}$ such that, for
all
$k \in \N$, for all
$k \leq n_1 <\ldots < n_k$, for all norm $1$ vectors $x_i$ in $G_{n_i}$, $i \leq k$, all
coefficients
$(\lambda_i)_{i \leq k}$, with $|\lambda_i| \leq 1$,
$$|f(\sum_{i=1}^k \lambda_i x_i)- \tilde{f}(\lambda_1,\ldots,\lambda_k)| \leq
\epsilon_k.$$

We recall that a basic sequence $(x_n)_{n \in \N}$ generates a spreading model
$(\tilde{x}_n)_{n \in \N}$ if
for any $\epsilon>0$, and $k \in \N$, there exists $N \in \N$ such that for all
$N<n_1<\ldots<n_k$, the sequences $(x_{n_i})_{i \leq k}$ and $(\tilde{x}_i)_{i
  \leq k}$ are $1+\epsilon$-equivalent. A spreading model is a basic sequence
which is necessarily $1$-subsymmetric.

The main application given in \cite{OR} for
 their result is about spreading models, and we derive from this the following lemma.

\begin{lemm} \label{FDD} Let $X$ be a Banach space with a Schauder
 (resp. unconditional) basis
$(e_n)_{n
\in \N}$. Let $p \in [1,+\infty]$ be such that $\ell_p$ is 
block (resp. disjointly) finitely representable in
$X$.  Then there 
exist a spreading model $(\tilde{y}_n)_{n \in \N}$ generated by a
block-sequence in $X$, a normalized
block-sequence (resp. sequence of disjointly supported blocks)
$(x_n)$ in
$X$, successive intervals $I_k$ forming a partition of $\N$ and some 
$C \geq 1$ such that:

\begin{itemize}
\item $|I_k|=k$ for all $k \in \N$,

\item for all $k \in \N$, $(x_n)_{n \in I_k}$ is
$C$-equivalent to the unit basis of $\ell_p^{k}$,

\item for any $k \in \N$, any $k<n_1<\ldots<n_k$, any normalized sequence
$(y_i)_{1 \leq i \leq k}$ with $(y_i) \in [x_n]_{n \in I_{n_i}}$, 
$\forall i \leq k$, the
sequence
$(y_i)_{1 \leq i \leq k}$ is $2$-equivalent to $(\tilde{y}_i)_{1 \leq i \leq k}$.

\end{itemize}

\end{lemm}

\pf Assume $\ell_p$ is block finitely representable in $X$.
We construct a block-subspace of $X$ of the form $\oplus_{n \in \N} F_n$, 
where each $F_n$ is a
block-subspace of dimension $n$ whose basis is $C$-equivalent to the 
basis of $\ell_p^n$ and the $F_n$'s are successive. We
apply the result of \cite{OR} to $\oplus_{n \in \N} F_n$ with the norm on $X$, which is a
Lipschitz map on $X$.

We pick a sequence $(\epsilon_k)$ of positive real numbers smaller than $1$
and decreasing to $0$. We obtain finite block-subspaces $G_k$ and a spreading model
$\tilde{y}_n$ such  that for any $k \in \N$,
any
$k<n_1<\ldots<n_k$, any normalized sequence
$(y_i)_{1 \leq i \leq k}$ with $(y_i) \in G_i$ is $1+\epsilon_k$-equivalent to
$\tilde{y}_n$.
We let $(x_n)_{n \in I_k}$ be the canonical basis of $G_k$ for all $k$ and we pass to a
subsequence to obtain the correct dimension $k$ for each $G_k$: $(x_n)_{n \in I_k}$ is
uniformly equivalent to the basis of $l_p^k$.

In the case when $\ell_p$ is disjointly finitely representable in $X$, we do the
same construction with the difference that each $F_n$ will have a
basis $C$-equivalent to the basis of $\ell_p^n$ which is disjointly supported on $X$, instead of successive. \pff

\begin{lemm} \label{lpsaturated} Let $X$ be a Banach space with an
unconditional basis. Assume $E_0$ is
not Borel reducible to permutative equivalence on $bb(X)$ (resp. on $dsb(X)$) and let 
$(f_n)_{n \in \N}$ be a Rosenthal
basic sequence in $X$ given by Lemma \ref{rosenthal}. Let $p \in 
[1,+\infty]$ be such that
$\ell_p$ is block-finitely representable in $[f_n]_{n \in \N}$ (resp. disjointly
finitely representable in $X$). Then 
$(f_n)_{n \in \N}$ is equivalent to
the unit basis of
$\ell_p$ (or $c_0$ if $p=+\infty$).
\end{lemm}

\pf Let $(f_n)$ be a Rosenthal basic sequence in $X$.  Let $p$ be 
such that $\ell_p$ is
block-finitely representable in $[f_n]_{n \in \N}$ (resp. disjointly finitely representable
in $X$). Let
$(e_n)$ be the canonical  basis of $\ell_p$ (or
$c_0$ if $p=+\infty$).  We need to prove that
$(f_n)$ is equivalent to $(e_n)$.

We note that any spreading model $(\tilde{y}_n)$ generated by a block-sequence in $X$ is
equivalent to $(f_n)$. Indeed, any block-sequence
 generating this spreading model has a 
subsequence equivalent to
$(f_n)$, so $(\tilde{y}_n)$ is equivalent to $(f_n)$. So
by Lemma \ref{FDD}, we find a block-sequence of $[f_n]_{n \in \N}$ (resp. sequence of
disjointly supported blocks of $X$) $(x_n)_{n \in \N}$, a constant 
$C \geq 1$
and associated intervals $(I_k)$ of length $k$ so that

\begin{itemize}
\item for all $k \in \N$, $(x_n)_{n \in I_k}$ is
$C$-equivalent to $(e_n)_{n \in I_k}$,

\item for any $k \in \N$, any $k<n_1<\ldots<n_{k}$, any normalized sequence
$(y_i)_{1 \leq i \leq k}$ with $(y_i) \in [x_n]_{n \in I_{n_i}}$, 
$\forall i \leq k$, the
sequence
$(y_i)_{1 \leq i \leq k}$ is $C$-equivalent to $(f_i)_{1 \leq i \leq k}$.
\end{itemize}

Passing to a subsequence of $(x_n)$,
 we may assume that 
$x_n<f_n^{\prime}<x_{n+1}$ for all $n \in
\N$, for some subsequence $(f_n^{\prime})$ of $(f_n)$
(resp. that $x_n$ and $f_p^{\prime}$ are disjointly
supported for all
$n,p$ in
$\N$). Applying Proposition
\ref{lppluslq} to
$(x_n)$ and
$(f_n^{\prime})$, and using the fact that $(f_n)$ is subsymmetric, we find a linear order
$\ext$ on $\N$ such that $(x_n) \leq^{C'} (f_n)^{\ext}$, for some  
constant $C'$. In
particular,
  for all $k \in \N$,
$$(x_n)_{n \in I_k} \leq^{C'} (f_n)^{\ext}_{n \in I_k}.$$ This
implies that
$$(e_n)_{n \leq k} \leq^{cCC'} (f_{\sigma(n)})_{n \leq k},$$ where $c$ is such that
$(f_n)$ is $c$-subsymmetric and $\sigma$ is a 
permutation
 on $\{1,\ldots,k\}$.
By symmetry of the basis $(e_n)$ and as $k$ was arbitrary, we deduce 
that $(f_n)$
$cCC'$-dominates $(e_n)$.

\

We now prove that $(f_n)$ is dominated by $(e_n)$, and to simplify 
the notation, we assume
$p<+\infty$; the case $p=+\infty$ is similar.  Assume on the contrary that
$(f_n)$ is not dominated by $(e_n)$. Then we may build by induction a 
rapidly converging
system $(\delta_k),(I_k)$ and some increasing sequence $K_k$ such 
that for all $k \in
\N$,

\begin{itemize}
\item (6) $K_k>2\sum_{j=1}^{k-1} |I_j|$ and $K_k \delta_k \geq k$,

\item (7)  there exists a sequence
$(\mu_i)_{i \in I_k}$ which satisfies
$\norm{\sum_{i \in I_k} \mu_i e_i} \leq 1$ and $\norm{\sum_{i \in I_k}
\mu_i f_i}
\geq K_k$.
\end{itemize}

We consider the previously defined sequence $(x_n)_{n \in \N}$ and, 
up to passing to a
subsequence of $(x_n)_{n \in \N}$ corresponding to the partition $(I_k)$, we 
may assume that for some subsequence $(f^{\prime}_n)$ of $(f_n)$,
$x_n<f_n^{\prime}<x_{n+1}$ for all
$n \in \N$ (resp. $x_n$ and $f_p^{\prime}$ are disjointly supported for all $n,p$ in
$\N$) , and that we have:

\begin{itemize}
\item for all $k \in \N$, $(x_n)_{n \in I_k}$ is
$C$-equivalent to $(e_n)_{n \in I_k}$.
\item for any $k \in \N$, any $k<n_1<\ldots<n_k$, any normalized sequence
$(y_i)_{1 \leq i \leq k}$ with $(y_i) \in [x_n]_{n \in I_{n_i}}$, 
the sequence
$(y_i)_{1 \leq i \leq k}$ is $C$-equivalent to $(f_i)_{1 \leq i 
\leq k}$.
\end{itemize}

  By Lemma \ref{lppluslqtechnical} applied to $(f_n^{\prime})$ and $(x_n)$ we 
may find  $D \geq 1$, an infinite subset $K$ of $\N$, and for all $k \in K$, a 
subset $J_k$ of $I_k$
with
$|I_k \setminus J_k| \leq
\sum_{j=0}^{k-1}|I_j|$ and distinct integers $(n_i)_{i \in J_k}$ such
that
$$(\delta_k f_i^{\prime}+x_i)_{i \in J_k} \sim^D (x_{n_i})_{i \in J_k}.$$
The end of our proof now divides in two cases.
For $k \in K$, let $A_k$ be the set of $n$'s
such that
$\{n_i, i \in J_k\}
\cap I_n \neq \emptyset$.

{\em First case:} we first assume that for any $m \in \N$, we may find $k \in K$ such that
the  set $A_k$  is of cardinal at least $m$.

Let $m \in \N$. For infinitely many $k$'s, we may find a set $L_k 
\subset J_k$ of
cardinal
$m$ such that
$\{n_i, i \in L_k\}$ meets $I_n$ for exactly $m$ values of $n$ which are 
strictly larger than $m$.
Then
$(x_{n_i})_{i \in L_k} \sim^C (f_{n_i})_{i \in L_k}.$

We deduce that
$$(f_{n_i})_{i \in L_k} \leq^{CD} (\delta_k f_i^{\prime} +x_i)_{i \in L_k}$$
so, as $L_k \subset I_k$, for all $(\lambda_i)_{i \in L_k}$,
$$\norm{\sum_{i \in L_k} \lambda_i f_{n_i}} \leq CD (\delta_k m 
\max_{i \in L_k}|\lambda_i|+
C(\sum_{i \in L_k}|\lambda_i|^p)^{1/p}),$$
and by symmetry of the basis of $\ell_p$ and $c$-subsymmetry of $(f_n)$, we 
deduce that for
any sequence $(\lambda_i)_{1 \leq i \leq m}$,

$$\norm{\sum_{i \leq m} \lambda_i f_{i}} \leq cCD(\delta_k m \max_{i \leq m}|\lambda_i|+C
(\sum_{i \leq m}|\lambda_i|^p)^{1/p}).$$
Letting $k$ tend to infinity and as $m$ was arbitrary, we obtain that
$(f_n)$ is $c C^2 D$-dominated by $(e_n)_{n \in \N}$.

{\em Second case:} we now assume that there exists some $m \in \N$ such that for all $k 
\in K$, the set
$A_k$ contains at most $m$ elements.

Then for any $k \in K$, all $(\lambda_i)_{i \in J_k}$,
$$\norm{\sum_{i \in J_k} \lambda_i
x_{n_i}}=
\norm{\sum_{n \in A_k} \sum_{i \in J_k, n_i \in I_n} \lambda_i x_{n_i}} \leq Cm
(\sum_{i \in J_k}|\lambda_i|^p)^{1/p}.$$

It follows that

$$\delta_k \norm{\sum_{i \in J_k}\lambda_i f_i^{\prime}} \leq CDm (\sum_{i \in
J_k}|\lambda_i|^p)^{1/p}.$$

Applying this to the coefficients $\mu_i$ given by (7), we obtain
$$\delta_k (K_k-\sum_{j=1}^{k-1}|I_j|) \leq cCDm,$$
where $c$ is such that $(f_n)$ is $c$-subsymmetric, so by (6),

$$k \leq \delta_k K_k \leq 2cCDm,$$ a contradiction. \pff

\begin{theo}\label{resultatblock} Let $X$ be a Banach space with a Schauder basis 
$(e_n)$. Assume
$E_0$ is not Borel reducible to permutative
equivalence on $bb(X)$. Then
 there exists $p \in [1,+\infty]$ such that every 
block-sequence
of
$X$ has a block-sequence which is equivalent to the canonical basis of
$\ell_p$ (or $c_0$ if $p=+\infty$).

\end{theo}

  \pf If $(e_n)$ is unconditional, Lemma \ref{rosenthal} applies, so there is a Rosenthal
basic sequence
 $(f_n)$ such that every normalized
block basis in $X$ has a subsequence equivalent to $(f_n)$. Let $p$ be such 
that $\ell_p$ is block
finitely representable in
$[f_n]_{n \in \N}$ ($p$ exists by Krivine's Theorem) . By Lemma \ref{lpsaturated}, $(f_n)$
is equivalent to the basis of
$\ell_p$ (or
$c_0$
 if $p=+\infty$).

In the general case, note that by \cite{R}, Proposition 21, every normalized
block-sequence in $X$ has a subsequence which is permutatively equivalent to its further
subsequences. In particular, $X$ contains no hereditarily indecomposable subspace (no
subspace of a H.I. space is isomorphic to a proper subspace \cite{GM}), and by Gowers'
dichotomy theorem,
$X$ is saturated with unconditional block-sequences.

By the unconditional case, we deduce that $X$ is saturated with spaces isomorphic to $c_0$
or $\ell_p$. Finally, if $X$ contains $\ell_p$ and $\ell_q$, for $p \neq q$, then as $\ell_p$ and $\ell_q$
are totally incomparable, $X$ contains a direct sum $\ell_p \oplus \ell_q$, and we may assume
that these copies are spanned by block-sequences $(x_n)$ and $(y_n)$ which alternate
(i.e. $\forall n \in \N, x_n<y_n<x_{n+1}$). By
Proposition
\ref{b2},
$E_0$ is Borel reducible to permutative equivalence on $bb_2(\ell_p \oplus \ell_q)$, so $E_0$
would be Borel reducible to $\sim^{perm}$ on $bb(X)$, a contradiction. The same proof
holds for $c_0$ and $\ell_p$. We deduce that there is a unique $p$ such that $X$ contains
copies of $\ell_p$ (or $c_0$ if $p=+\infty$).
  \pff

\

A Banach space $X$ with an unconditional basis is said to satisfy an {\em upper $p$
estimate} if there
exists $C \geq 1$ such that for any disjointly supported vectors
$x_1,\ldots,x_n$, $\norm{\sum_{i=1}^n x_i} \leq
 C(\sum_{i=1}^n \norm{x_i}^p)^{1/p}$ (or $\norm{\sum_{i=1} x_i} \leq C\max_{i \leq n}\norm{x_i}$ if $p=+\infty$). By an easy uniform boundedness
argument,
 this
is equivalent to saying that for any normalized disjointly supported sequence
$(x_n)_{n \in \N}$ on $X$, $(x_n)$ is dominated by the canonical basis of
$\ell_p$ (or $c_0$ if $p=+\infty$).

\begin{theo}\label{resultat} Let $X$ be a Banach space with an unconditional
 basis 
$(e_n)$.
\begin{itemize}
\item
 Assume
$E_0$ is not Borel reducible to permutative
equivalence on $bb(X)$. Then there exists $p \in [1,+\infty]$ such that
every normalized block-sequence
of
$X$ has a subsequence which is equivalent to the canonical basis of $\ell_p$ (or $c_0$ if
$p=+\infty$).

\item
 Assume
$E_0$ is not Borel reducible to permutative
equivalence on $dsb(X)$. Then there is a unique $p \in [1,+\infty]$ such that
$\ell_p$ is disjointly finitely representable in $X$. If $p=+\infty$ 
then $(e_n)$ is
equivalent to the unit vector basis of $c_0$. If $p<+\infty$ then $X$ 
satisfies an upper
$p$-estimate and
every normalized block-sequence
of
$X$ has a subsequence which is equivalent to the canonical basis of $\ell_p$.
\end{itemize}
\end{theo}

\pf The $bb(X)$ case is proved at the beginning of the proof of Theorem \ref{resultatblock}.
Assume now that $E_0$ is not Borel reducible to permutative
equivalence on $dsb(X)$.
By Lemma \ref{rosenthal}, there
is a Rosenthal basic sequence
 $(f_n)$, necessarily unique up to equivalence, such that every normalized
block basis in $X$ has a subsequence equivalent to $(f_n)$. Let $p$ be such 
that $\ell_p$ is disjointly 
finitely representable in
$X$. By Lemma \ref{lpsaturated}, $(f_n)$ is equivalent to
the basis of
$\ell_p$ (or
$c_0$
 if $p=+\infty$), so $p$ is unique. It remains to show
that
$(e_n)$ satisfies
 an upper $p$-estimate, which implies that $(e_n)$ is equivalent to the basis
 of $c_0$ if $p=+\infty$.

For any $(x_n) \in dsb(X)$, we may find a normalized
 sequence $(v_n) \sim (f_n)$ which is
disjointly supported from $(x_{2n})$. As $E_0$ is not Borel reducible to
permutative
equivalence on $bb_2([x_{2n}]\oplus[v_n])$, we deduce from
Proposition \ref{lppluslq} that $(x_{2n}) \leq (v_n^{\ext})$ for some linear order
$\ext$ on $\N$. As $(v_n)$ is symmetric it follows that
$(x_{2n}) \leq (v_n)$, that is for some $C$ and any $(\lambda_n) \in c_{00}$,
$$\norm{\sum_{n \in \N} \lambda_{2n}x_{2n}} \leq C(\sum_{n \in
  \N}|\lambda_{2n}|^p)^{1/p},$$
if $p<+\infty$, or
$$\norm{\sum_{n \in \N} \lambda_{2n}x_{2n}} \leq C\max_{n \in
  \N}|\lambda_{2n}|,$$
 if $p=+\infty$.
We obtain a similar estimate for $(x_{2n+1})$ and  deduce that $(x_n)$ is
  dominated by the unit vector basis of $\ell_p$ (or $c_0$ if $p=+\infty$), and so
finally $(e_n)$ satisfies an upper $p$-estimate.\pff

\

Note that from this theorem, we may deduce that $E_0$ is Borel reducible to $\sim^{perm}$
on $bb(S)$, where $S$ is Schlumprecht's space \cite{S}. It is however still unknown if $S$
is ergodic.

\section{Permutative equivalence between uncondi\-tional ba\-sic 
se\-quences in $X$ and in $X^*$.}

We obtain a complete dichotomy result by looking at the disjointly supported sequences of
the dual $X^*$
of $X$, when $X^*$ has a basis. Compare this theorem with Conjecture \ref{conjectureperm}.

\begin{theo} \label{dsb} Let $X$ be a Banach space with a shrinking normalized
unconditional 
  basis
$(e_n)$. Then either $(e_n)$ is equivalent to the canonical basis of $c_0$ or 
some $\ell_p, 1<p<+\infty$,
or
$E_0$ is Borel reducible to permutative
equivalence on $dsb(X)$,
or on $dsb(X^*)$.
\end{theo}

\pf Assume $E_0$ is Borel reducible to permutative equivalence 
neither on $dsb(X)$ nor on
$dsb(X^*)$. By Theorem \ref{resultat},  there exists $Y=c_0$ or $\ell_p$ for some 
$1<p<+\infty$ such that
every normalized block-sequence
  of $X$ has a subsequence
equivalent to the canonical basis of $Y$, and
we may assume that $1<p<+\infty$ and that $X$ satisfies 
an upper $p$-estimate.

   Some subsequence of $(e_n)$ is 
equivalent to the basis
of $\ell_p$, so its dual basis identified with a subsequence of 
$(e_n^*)$ is equivalent to
the basis of $\ell_{p'}$ (where $\frac{1}{p}+\frac{1}{p'}=1$). Thus by 
 Theorem \ref{resultat} applied for $X^*$,
$X^*$ satisfies an upper $p'$-estimate.
So $(e_n^*)_{n \in \N}$ is dominated by the unit vector basis of $\ell_{p'}$. It follows that
$(e_n)_{n \in \N}$ dominates the unit vector basis of $\ell_p$.
Finally $(e_n)_{n \in \N}$ is equivalent to the unit vector basis of $\ell_p$.\pff

\

We also deduce the following dichotomy result about the number of non permutatively
equivalent
sequences spanning subspaces, or quotients, of a Banach space with an
unconditional basis which is not
isomorphic to a Hilbert space. Note that by uniqueness of the unconditional
basis of $\ell_2$, any normalized unconditional basis of a subspace, or a quotient, of
$\ell_2$ must be (permutatively) equivalent to the canonical basis of $\ell_2$.

\begin{theo}\label{continuum} Let $X$ be  a Banach space with an unconditional basis. Then
  either $X$  is  isomorphic
to $\ell_2$, or
$X$ contains
$2^{\omega}$ subspaces, or $2^{\omega}$ quotients, spanned by normalized 
 unconditional bases which
are mutually non permutatively equivalent.
\end{theo}

\pf Assume $X$ is not isomorphic to $\ell_2$. If $X$ contains $c_0$ or 
$\ell_1$, then we are done,
since by
\cite{FG}, there is a Borel reduction of $E_0$ to permutative 
equivalence between the
canonical unconditional bases of some subspaces of
$c_0$ (resp
$\ell_1$). So by the classical result of James (see \cite{LT}),
 we may assume $X$ is reflexive.

We may assume the basis of $X$ is normalized and we apply Theorem \ref{dsb}. If $E_0$ is
Borel reducible to permutative  equivalence on $dsb(X)$,
then we obtain the desired result with subspaces of $X$. If $E_0$ is 
Borel reducible to
permutative equivalence on $dsb(X^*)$, let $f:2^{\omega} \rightarrow 
dsb(X^*)$ be the Borel
reduction. We note that the bases
$f(\alpha)$ and $f(\beta)$ are permutatively equivalent if and only 
if the dual bases
$f(\alpha)^*$ and $f(\beta)^*$ are permutatively equivalent; and for 
$\alpha \in
2^{\omega}$,
the dual basis $f(\alpha)^*$ is an unconditional basis of some 
quotient of $X$. We thus
obtain continuum many non permutatively equivalent normalized unconditional
 bases of quotients 
of $X$ in the family
$f(\alpha)^*, \alpha \in 2^{\omega}$.

Finally if the basis of $X$ is equivalent to the canonical basis of
 some $\ell_p, 1 <p 
<+\infty$, with $p<2$,
\cite{FG} gives an explicit construction of $2^{\omega}$ subspaces of $X$ with normalized
unconditional bases which are mutually non permutatively equivalent (see the
 proof of Proposition \ref{bossard}; in fact we even obtain a reduction of
 $E_0$ to permutative equivalence between such
 unconditional bases in that case). 
If $p>2$, then we use
duality to deduce the existence of $2^{\omega}$ quotients of $X$ with normalized
unconditional bases
which are mutually non permutatively equivalent.\pff

\

The reader should compare this result with Conjecture \ref{l2}, noting that the proof of
Theorem
\ref{continuum} actually gives a reduction of $E_0$ to permutative equivalence on
an appropriate space of basic sequences spanning subspaces or quotients of $X$, when $X$ is
not isomorphic to $\ell_2$.

To conclude, let us mention two results with some similarity with Theorem \ref{continuum},
by the use their hypotheses make of both subspaces and duals (resp. quotients). By P.
Mankiewicz and N. Tomczak-Jaegermann, if every subspace of every quotient of $\ell_2(X)$
has a Schauder basis, then the Banach space $X$ must be isomorphic to $\ell_2$ \cite{MTJ}.
By V. Mascioni, if
$\ell_2(X)$ is locally self-dual (i.e. finite dimensional subspaces are uniformly
isomorphic to their duals), then $X$ must also be isomorphic to $\ell_2$ \cite{Ma}.

\

\

 Equipe d'Analyse Fonctionnelle,

Universit\'e Paris 6,

Bo\^\i te 186, 4, Place Jussieu,

 75252, Paris Cedex 05,

 France.

\

 E-mail: ferenczi@ccr.jussieu.fr.

Phone: (33)-1-44-27-54-34

Fax :(33)-1-44-27-25-55


\begin{thebibliography}{Wwww}




\bibitem{S} G. Androulakis and T. Schlumprecht, {\em The
Banach space
$S$ is complementably minimal and subsequentially prime},  Studia Math.  {\bf 156}  (2003), 
3, 227--242.





\bibitem{B} B. Bossard, {\em A coding of separable Banach spaces.
Analytic and coanalytic families of Banach spaces.} Fund. Math. {\bf
172}
 (2002), no. 2, 117--152.





\bibitem{BCLT} J. Bourgain, P. Casazza, J. Lindenstrauss and L.
Tzafriri, {\em Banach spaces
with a unique unconditional basis, up to permutation}, Memoirs of the
A.M.S. {\bf 54} (1985), No 322.



\bibitem{CK1} P. Casazza and N.J. Kalton, {\em Uniqueness of unconditional bases in
    Banach spaces}, Israel J. Math. {\bf 103} (1998), 141--176.



\bibitem{CK2} P. Casazza and N.J. Kalton, {\em Uniqueness of unconditional bases in
    $c_0$-products}, Studia Math {\bf 133} (1999), 3, 275--294.


\bibitem{F}
V. Ferenczi,
{\em Minimal subspaces and isomorphically homogeneous sequences in a
Banach space},
Israel Journal of Mathematics, to appear.

\bibitem{F2}
V. Ferenczi,
{\em Topological 0-1 laws for
subspaces of a Banach
space with a Schauder basis}, preprint.

\bibitem{FG}
V. Ferenczi and E. M. Galego, {\em Some equivalence relations which are Borel
reducible to 
 isomorphism  between separable
Banach spaces}, Israel Journal of Mathematics, to appear.

\bibitem{FPR} V. Ferenczi, A.M. Pelczar, and C. Rosendal, {\em On a
question of H. P. Rosenthal
concerning a characterization of $c_0$ and $\ell_p$}, Bull. London Math.
Soc. {\bf 36} (2004), 3,
396--406.

  

\bibitem{FR1} V. Ferenczi and C. Rosendal,
 {\em On the number of non-isomorphic subspaces of a Banach space},
 Studia Math.  {\bf 168}  (2005),  no. 3, 203--216.












\bibitem{FR2} V. Ferenczi and C. Rosendal, {\it Ergodic Banach
spaces}, Advances in Mathematics {\bf 195} (2005), 1, 259--282.







\bibitem{FS} H. Friedman and   L. Stanley, {\em
A Borel reducibility theory for classes of countable structures.}
J. Symbolic Logic {\bf 54} (1989), no. 3, 894--914.







\bibitem{G} W.T. Gowers, {\em An infinite Ramsey theorem and some
Banach-space dichotomies.}, Ann. of
Math. (2) {\bf 156} (2002), no. 3, 797--833.



\bibitem{G2} W.T. Gowers, {\em A solution to the Schroeder-Bernstein problem
    for Banach spaces}, Bull. London Math. Soc. {\bf 28} (1996), 297--304.



\bibitem{GM} W.T. Gowers and B. Maurey, {\em The unconditional basic
sequence problem}, J. Amer.
Math. Soc. {\bf 6} (1993),4, 851--874.



\bibitem{GM2} W.T. Gowers and B. Maurey, {\em Banach spaces with small spaces
    of operators}, Math. Ann. {\bf 307} (1997), 543--568.



\bibitem{HKL} L. A. Harrington, A. S. Kechris and   A. Louveau, {\em
A Glimm-Effros dichotomy for Borel equivalence relations.}
J. Amer. Math. Soc. {\bf 3} (1990), no. 4, 903--928.



\bibitem{K} N. Kalton, {\em Lattice structures on Banach spaces},
  Mem. Amer. Math. Soc. {\bf 103},
  1993, no 493.





\bibitem{Ke}  A. S. Kechris, {\it Classical descriptive set theory.}
Graduate Texts in Mathematics,
{\bf 156}, Springer-Verlag, New York, 1995.







   \bibitem{KT}  R. A. Komorowski and N. Tomczak-Jaegermann, {\em
Banach spaces without local
unconditional structure}. Israel J. Math. {\bf 89}
   (1995), no. 1-3, 205--226.   {\em Erratum to ``Banach spaces without
local unconditional structure''}. Israel J. Math. {\bf 105} (1998),
85--92.









\bibitem{L} H. Lemberg, {\em Nouvelle d\'emonstration d'un 
th\'eor\`eme de J.L. Krivine sur
la finie repr\'esentation de $\ell_p$ dans un espace de Banach}, Israel J. Math.
{\bf 39} (1981), 341--348.



   \bibitem{LT} J. Lindenstrauss and L. Tzafriri, {\it Classical Banach
spaces}, Springer-Verlag New
York Heidelberg Berlin. 1979.

\bibitem{MTJ} P. Mankiewicz and N. Tomczak-Jagermann, {\em Schauder bases in quotients
of subspaces of
$\ell_2(X)$}, Amer. J. Math. {\bf 116} (1994), no. 6, 1341--1363.


\bibitem{Ma} V. Mascioni, {\em On Banach spaces isomorphic to their duals},  Houston J.
Math.  {\bf 19}  (1993),  no. 1, 27--38.


\bibitem{M} B.S. Mityagin, {\em Equivalence of bases in Hilbert scales (in
    Russian)}, Studia Math. {\bf 37} (1970), 111--137.




\bibitem{OR} E. Odell, H.P. Rosenthal and  Th. Schlumprecht, {\em On weakly null
FDDs in Banach spaces}, Israel J. Math {\bf 84} (1993),3, 333--351.

\bibitem{P} A. Pe\l czy\'nski, {\em Universal bases}, Studia Math. {\bf 32} (1969),
247--268.

\bibitem{Rthese} C. Rosendal, {\em Etude descriptive de l'isomorphisme dans la classe
des espaces de Banach}, Th\`ese de Doctorat de l'Universit\'e Paris 6 (2003).

\bibitem{R}  C. Rosendal, {\em Incomparable, non isomorphic and minimal Banach
      spaces}, Fund. Math.  {\bf 183}  (2004) 3, 253--274.


\bibitem{R1} C. Rosendal, {\em Cofinal families of Borel equivalence relations and
quasiorders}, J. Symbolic Logic, to appear.





\bibitem{W} P. Wojtaszczyk, {\em Uniqueness of unconditional bases in
    quasi-Banach spaces with applications to Hardy spaces, II}, Israel J. Math
{\bf 97} (1997), 253--280.



\bibitem{Wo} M.  Wojtowicz, {\em On Cantor-Bernstein type theorems in Riesz
    spaces}, Indag. Math. {\bf 91} (1998), 93--100.



\end{thebibliography}
\end{document}